\documentclass[12pt]{article}

\usepackage{amsmath}
\usepackage{amsfonts}
\usepackage{amsthm}
\usepackage{graphicx}
\usepackage{caption}
\usepackage{subcaption}
\usepackage{enumerate}
\usepackage{hyperref}
\usepackage[parfill]{parskip} 

\usepackage{fancyhdr}
\numberwithin{equation}{section}

\begin{document}

\title{Explicit (Polynomial!) Expressions for the Expectation, Variance and Higher Moments of the Size of a $(2n+1,2n+3)$-core partition with Distinct Parts }
\author{Anthony Zaleski \and Doron Zeilberger}

\maketitle

\begin{abstract}
Armin Straub's beautiful article (\url{https://arxiv.org/abs/1601.07161})
concludes with two intriguing conjectures about the number, and maximal size, of $(2n+1,2n+3)$-core partitions 
with distinct parts. These were proved by  ingenious, but complicated, arguments
by Sherry H.F. Yan, Guizhi Qin, Zemin Jin, and Robin D.P. Zhou (\url{https://arxiv.org/abs/1604.03729}). 
In the present article, we first comment that these results can be proved faster by ``experimental mathematics'' methods, 
that are easily rigorizable.
We then develop relatively efficient, symbolic-computational, algorithms, based on non-linear functional recurrences,
to generate what we call the Straub polynomials, where $S_n(q)$ is  the generating function, according to size,
of the  set of $(2n+1,2n+3)$-core partitions with distinct parts, and compute the first $21$ of them.
These are used to deduce explicit expressions, as polynomials in $n$, for the mean, variance, and the third through the 
seventh moments (about the mean) of the random variable ``size'' defined on $(2n+1,2n+3)$-core partitions with 
distinct parts. 
In particular we show that this random variable is not asymptotically normal, and the limit of the coefficient of variation is 
$\sqrt{14010}/150 = 0.789092305...$,
the scaled-limit of the third moment (skewness) is  $(396793/390815488)\cdot \sqrt{467 \cdot 7680}= 1.92278748..$, 
and that the scaled-limit of the 4th-moment (kurtosis) is 
$145309380/16792853= 8.6530490...$. 
We are offering to donate one hundred dollars to the OEIS foundation in honor of the first to identify the limiting distribution. 
\end{abstract}
\leavevmode
\\
\\
\\
\section*{Supporting Maple packages and output}

All the results in this article were obtained by the use of the Maple packages

$\bullet$ {\tt http://www.math.rutgers.edu/\~{}zeilberg/tokhniot/Armin.txt}  ,

$\bullet$ {\tt http://www.math.rutgers.edu/\~{}zeilberg/tokhniot/core.txt}  ,

whose output files, along with links to diagrams, are available from the {\it front} of this article

{\tt http://www.math.rutgers.edu/\~{}zeilberg/mamarim/mamarimhtml/armin.html}.

 \section{$(s,t)$-core partitions and drew Armstrong's ex-conjecture}

Recall that a {\it partition} is a non-increasing sequence of positive integers $\lambda=(\lambda_1, \dots, \lambda_k)$ with $k \geq 0$,
called its {\it number of parts}; $n:=\lambda_1 + \dots + \lambda_k$ is called its {\it size},
and we say that $\lambda$ is a {\it partition of $n$}.
Also recall that the {\it Ferrers diagram} 
(or equivalently, using empty squares rather than dots, {\it Young diagram}) 
of a partition $\lambda$ is obtained by placing, in a {\it left-justified} way,
$\lambda_i$ dots at the $i$-th row. For example, the Ferrers diagram of the partition $(5,4,2,1,1)$ is
$$
\begin{array}{c c c c c}
* & * & *  &*  & * \\
* & * & *  &*  &  \\
* & * &   &  &  \\
* &  &   &  &  \\
* &  &   &  &  \\
\end{array}  .
$$

Recall also that the {\it hook length} of a dot $(i,j)$ in the Ferrers diagram, $1\leq j \leq \lambda_i$, is the
number of dots to its right (in the same row) plus the number of dots below it (in the same column) plus one
(for itself), in other words $\lambda_i -i+\lambda'_j-j+1$, where $\lambda'$ is the {\it conjugate partition}, 
obtained by reversing the roles of rows and columns. (For example
if $\lambda=(5,4,2,1,1)$ as above, then $\lambda'=(5,3,2,2,1)$).

Here is a table of hook-lengths of the above partition, $(5,4,2,1,1)$:

$$
\begin{array}{c c c c c}
9 & 6 & 4  & 3  & 1 \cr
7 & 4 & 2  &1  &  \cr
4 & 1 &   &  &  \cr
2 &  &   &  &  \cr
1 &  &   &  &  \cr
\end{array}   .
$$
It follows that its set of hook-lengths is $\{1,2,3,4,6,7,9\}$. A partition is called an $s$-core if none of its hook-lengths is $s$.
For example, the above partition, $(5,4,2,1,1)$, is a $5$-core, and an $i$-core for all $i \geq 10$.

A partition is a {\it simultaneous} $(s,t)$-core partition if it avoids both $s$ and $t$. For example the
above   partition, $(5,4,2,1,1)$, is a $(5,11)$-core partition (and a $(5,12)$-core partition, and a $(100,103)$-core partition etc.).

For a lucid and engaging account, see [AHJ].

As mentioned in [AHJ], Jaclyn Anderson ([A]) very elegantly proved the following.

{\bf Theorem ([A]):} If $s$ and $t$ are relatively prime positive integers, then there are
{\bf exactly}
$$
\frac{(s+t-1)!}{s!t!} ,
$$
$(s,t)$-core partitions.

For example, here are the $(3+5-1)!/(3!5!)=7$  $(3,5)$-core partitions:
$$
\{ empty , 1, 2, 11, 31, 211, 4211 \}  .
$$

Drew Armstrong ([AHJ], conjecture 2.6) conjectured, what is now the following theorem.

{\bf Theorem ([J])} : The {\it average size} of an $(s,t)$-core partition is given by
the nice polynomial
$$
\frac{(s-1)(t-1)(s+t+1)}{24}  .
$$

For example, the (respective) sizes of the above-mentioned $(3,5)$-core partitions are
$$
0,1,2,2,4,4,8 ,
$$
hence the average size is
$$
\frac{0+1+2+2+4+4+8}{7}=\frac{21}{7}=3 ,
$$
and this agrees with Armstrong's conjecture, since
$$
\frac{(3-1)(5-1)(3+5+1)}{24} \, = \, 3  .
$$

Armstrong's conjecture was  proved by Paul Johnson ([J]) using a very complicated (but ingenious!) 
argument (that does much more).
Shortly after, and almost {\it simultaneously} (no pun intended) it was re-proved by Victor Wang [Wan],
using another ingenious (and even more complicated) argument, that also does much more, in particular,
proving an intriguing conjecture of Tewodros Amdeberhan and Emily (formerly Leven) Sergel ([AmL]).
Prior to the full proofs by  Johnson and Wang, Richard Stanley and Fabrizio Zanello [StaZ] came up with
a nice (but rather {\it ad hoc}) proof of the important special case of $(s,s+1)$-core partitions.
Explicit expressions for the {\it variance} and the third moment were found by Marko Thiel and  Nathan Williams ([TW]).

Ekhad and Zeilberger ([EZ]) went far beyond, and derived explicit expressions for the first 6 moments for the general
$(s,t)$-core partitions, and the first 9 moments for the case $(s,s+1)$, and used them to find
the scaled limits up to the ninth, that strongly suggest that the limiting distribution is
the continuous random variable
$$
\sum_{k=1}^{\infty} \frac{z_k^2 + \tilde{z}_k^2}{4 \pi^2 k^2}  ,
$$
where $z_k$ and $\tilde{z}_k$ are jointly independent sequences of independent standard normal random variables.

\section{Simultaneous core partitions into distinct parts}

Tewodros Amdeberhan ([Am]) initiated the study of simultaneous core partitions
with {\it distinct parts}, and conjectured that the number of $(s,s+1)$-core partitions with distinct parts is given by
the Fibonacci number $F_{s+1}$. This was proved by Armin Straub ([Str1]) and Huan  Xiong ([X]). Xion also
proved a conjectured expression of Amdeberhan for the expected size, in terms of a double sum involving
Fibonacci numbers. A more explicit expression was derived by the first-named author [Za], 
who also derived, assisted by his computer, {\bf explicit expressions} (as rational functions in $F_s, F_{s+1}$, and $s$) for the first $16$ moments.
He then deduced  that the scaled moments tend to the moments
of the standard normal distribution, giving strong evidence (that could be turned into a fully rigorous proof,
using the method of [Ze2])
that the random variable `size' defined over {\it distinct} $(s,s+1)$-core partitions is {\it asymptotically normal}. 

This is surprising, since, as already mentioned above,  it was
shown in [EZ] that when defined over all (not necessarily distinct) partitions, the random variable `size' is {\bf not} asymptotically normal.

At the end of his beautiful paper, [Str1], (where, among many other things, 
the author describes a beautiful
new elegant partition identity  between Odd and Distinct integer partitions which preserves the perimeter, that should have been found by Euler
(but had to wait for Straub)) Armin Straub conjectured two intriguing enumeration results.

{\bf Theorem 0} (conjectured in [Str1], first proved in [YQJZ]) The {\it number} of $(2n+1,2n+3)$-core partitions with distinct parts equals $4^n$.

{\bf Theorem 0'}: (conjectured in [Str1], first proved in [YQJZ])
The largest size of a $(2n+1,2n+3)$-core partition with distinct parts
is $\frac{1}{24}\, \left( 5\,n+11 \right) n \left( n+2 \right)  \left( n+1 \right)$.

The  proofs in [YQJZ] use
ingenious, but rather complicated, combinatorial arguments. We will, in this article, give new, much simpler,
`experimental-mathematical' proofs, that can be easily made rigorous. But our main purpose is to
establish {\it explicit} expressions for the expectation, variance, and all the moments up to the seventh.
With more computing power, it should be possible to go beyond. We then go on and use these explicit (polynomial) expressions in order to
find the limits of the scaled moments, giving exact values for the first seven moments of the limiting (scaled) probability distribution
of the random variable `size' over $(2n+1,2n+3)$-core partitions with distinct parts (as $n \rightarrow \infty$), and one of us (DZ) is pledging
\$100 to the OEIS foundation for identifying that limiting (continuous) probability distribution.

\subsection{Explicit expressions for the first seven moments}

{\bf Theorem 1}: The average size of a $(2n+1,2n+3)$-core partition with distinct parts is
$$
\frac{1}{32} ( 10\,{n}^{3}+27\,{n}^{2}+19\,n )  .
$$

Note that the corresponding average taken over {\it all}  partitions, according to Armstrong's ex-conjecture, is
$\frac{1}{6}n(n+1)(2n+5)=\frac{1}{3}n^3 +O(n^2)$, while, according to Theorem 1, our average (i.e. for the distinct case)
is $\frac{5}{16} n^3 +O(n^2)$, so it is a bit less.

{\bf Theorem 2}: The variance of the random variable `size' defined on the set of $(2n+1,2n+3)$-core partitions with distinct parts is
$$
\frac{1}{15360} (934\,{n}^{6}+4687\,{n}^{5}+9700\,{n}^{4}+10505\,{n}^{3}+6256\,{n}^{2}+1518\,n)  .
$$

Note that according to [EZ], the corresponding variance, taken over {\it all} partitions is
$$
{\frac {1}{720}}\, \left( 2\,n+1 \right)  \left( 2\,n+3 \right)  \left( 2\,n+2 \right) n \left( 4\,n+5 \right)  \left( 4\,n+4 \right) 
$$
which is  $\frac{8}{45} n^6 + O(n^5)= 0.1777777778 n^6 +O(n^5)$, while for our case, according to Theorem 2, it is
$\frac{467}{7680}n^6 +O(n^5)=0.06080729167 n^6 +O(n^5)$.

{\bf Theorem 3}: The  third moment (about the mean) of the random variable `size' defined on $(2n+1,2n+3)$-core partitions with distinct parts is
$$
\frac{1}{27525120} \cdot (
793586\,{n}^{9}+4945025\,{n}^{8}+12775144\,{n}^{7}+17215282\,{n}^{6}+11839450\,{n}^{5}$$
$$
+1535905\,{n}^{4}
-4756804\,{n}^{3}-4342612\,{n}^{2}-1297776\,n)  .
$$

{\bf Theorem 4}: The  fourth moment (about the mean) of the random variable `size' defined on $(2n+1,2n+3)$-core partitions with distinct parts is
$$
\frac{1}{54499737600} \cdot
(
1743712560\,{n}^{12}+13490284234\,{n}^{11}+45408125279\,{n}^{10}
$$
$$
+87568584895\,{n}^{9}+109173019890\,{n}^{8}+97494786972\,{n}^{7}+68082466947\,{n}^{6}
$$
$$
+34594762895\,{n}^{5}+8734303600\,{n}^{4}+3269131844\,{n}^{3}+7648567524\,{n}^{2}+4135638960\,n)  .
$$

{\bf Theorem 5}: The  fifth moment (about the mean) of the random variable `size' defined on $(2n+1,2n+3)$-core partitions with distinct parts is
$$
\frac{1}{108825076039680} \cdot
n \left( n+1 \right)  ( 4115597238066\,{n}^{13}+30331407775461\,{n}^{12}$$
$$
+93240357590320\,{n}^{11}+153901186416765\,{n}^{10}
+154511084293844\,{n}^{9}$$
$$
+126787455814599\,{n}^{8}+115227024155664\,{n}^{7}+42586120680111\,{n}^{6}
$$
$$
-95604599727502\,{n}^{5}
-105409116317640\,{n}^{4}+43165327777096\,{n}^{3}
$$
$$
+
91113907956144\,{n}^{2}-30975685518528\,n-65049004454400)  .
$$

{\bf Theorem 6}: The  sixth moment (about the mean) of the random variable `size' defined on $(2n+1,2n+3)$-core partitions with distinct parts is
$$
\frac{1}{8288117791182028800} \cdot
$$
$$
(
459077029253573970\,{n}^{18}+3986958940758529155\,{n}^{17}+14588638597341766281\,{n}^{16}
$$
$$
+29315654117562943844\,{n}^{15}+38855616058049391120\,{n}^{14}+ 52048632801161949890\,{n}^{13}
$$
$$
+87053992212835094382\,{n}^{12}+102228197171521441748\,{n}^{11}+24538654588404043230\,{n}^{10}
$$
$$
-81063397918244586845\,{n}^{9}-37681424022539337807\,{n}^{8}+128753068232342353072\,{n}^{7}
$$
$$
+136357236921377110920\,{n}^{6}-109095423240535042640\,{n}^{5}-264555566724556223856\,{n}^{4}
$$
$$
-62480060539123323264\,{n}^{3}+164786511770490504960\,{n}^{2}+100625844884387235840\,n
) \,\, .
$$


{\bf Theorem 7}: The  seventh moment (about the mean) of the random variable `size' defined on $(2n+1,2n+3)$-core partitions with distinct parts is
$$
\frac{n(n+1)}{
 2 ^{40} \cdot 3^5 \cdot 5^2 \cdot 7 \cdot 11 \cdot 13 \cdot 17 \cdot 19} \cdot
$$
$$
(203253344355858784830\,{n}^{19}+1525941518277673062635\,{n}^{18}
$$
$$
+4376090780890032310694\,{n}^{17}+5920532244827036954724\,{n}^{16}
$$
$$
+7108181147332994381598\,{n}^{15}+22516614862619041657440\,{n}^{14}
$$
$$
+47737754432542468750710\,{n}^{13}+21431538183386052191306\,{n}^{12}
$$
$$
-77127349790945221221652\,{n}^{11}-98788608530944679782107\,{n}^{10}
$$
$$
+91468628175188699900748\,{n}^{9}+276198594921821905993026\,{n}^{8}$$
$$
+164310592679893652073504\,{n}^{4} +1420837514400804031281984\,{n}^{3}
$$
$$
+53152679358583919475360\,{n}^{7}-516374679437475960870016\,{n}^{6}
$$
$$
-696941224296942655687312\,{n}^{5}
+1109985197630308975715328\,{n}^{2}
$$
$$
-745951061503715454673920\,n-1026387551269849288826880)  .
$$

\subsection{Corollaries}

{\bf 1.} The limit of the ``coefficient of variation'' (the quotient of the standard deviation to the mean),
as $n \rightarrow \infty$, is ${\frac {1}{150}}\,\sqrt {14010}=0.7890923055426827989 \dots$.
In particular, since that limit is {\bf not} zero,
unlike  $(k,k+1)$-core partitions with distinct parts discussed in [Za], there is {\bf no} concentration about the mean.

{\bf 2.} The limit of the {\it skewness}, as $n \rightarrow \infty$, is ${\frac {396793}{390815488}}\,\sqrt {467}\sqrt {7680}$ \\ $=1.922787480888358667\dots$

{\bf 3.} The limit of the {\it kurtosis}, as $n \rightarrow \infty$, is ${\frac {145309380}{16792853}}= 8.6530490084085\dots$

{\bf 4.} The limit of the scaled fifth moment ($\alpha_5$), as $n \rightarrow \infty$, is 
${\frac {3429664365055}{156594294624768}}\,\sqrt {467}\sqrt {7680}= 41.4777067204457\dots$

{\bf 5.} The limit of the scaled sixth moment ($\alpha_6$), as $n \rightarrow \infty$, is 
${\frac {382564191044644975}{1552893421695616}}=246.35572905\dots$.

{\bf 6.} The limit of the scaled seventh moment ($\alpha_7$), as $n \rightarrow \infty$, is 
${\frac {56459262321071884675}{62988906654652346368}}\,\sqrt {467}\sqrt {7680}=697.5015509357\dots$

\section{Proving the theorems}
We now explain the methods used to obtain the results in the previous section.
\subsection{A   New (``Experimental Math'') proof of Armin Straub's Ex-Conjecture that the number of $(2n+1,2n+3)$-core partitions
with distinct parts equals $4^n$}

The way Jaclyn Anderson proved  her celebrated theorem ([An]) that 
if $gcd(s,t)=1$, then the  number of  $(s,t)$-core partitions equals $(s+t-1)!/(s!t!)$ was by defining a bijection with the
set of {\bf order ideals} of the poset
$$
P_{s,t} := {\bf N} \backslash (s{\bf N} + t{\bf N} ) ,
$$
where ${\bf N}=\{ 0,1,2,3, \dots, \}$ is the set of non-negative integers, and the partial-order relation $c \, \leq_P \, d$ 
holds whenever $d-c$ can be expressed as $\alpha s + \beta t$ for some $\alpha,\beta \in {\bf N}$.

The set of order ideals of $P_{s,t}$, in turn, is in bijection with the set of {\it lattice paths} in the two-dimensional
square lattice, from $(0,0)$ to $(s,t)$ lying above the line $sy-tx=0$. 
This correspondence is used in the Maple package {\tt core.txt}, and was used in [Za], but
for our present purposes it is more efficient to use order ideals.

Recall that an order ideal $I$, in a poset $P$, is a set of vertices of $P$ such that if $c \in I$ then
all elements, $d$,  such that $d \,  \leq_P \, c$ also belong to I. Equivalently, if $d$ does {\bf not} belong to $I$, then
all vertices $c$ `above' it (i.e. such that $c \, \geq_P \, d$), also do {\bf not} belong to $I$.

Let $s(n)$ be the number of order ideals of the lattice $P_{2n+1,2n+3}$ with no  consecutive labels. 
Recall that, thanks to Jaclyn Anderson, this is the number of $(2n+1,2n+3)$-core partitions with distinct parts, our
{\it object of desire}.

Let's try and find an algorithm to compute the sequence
$\{s(n)\}$ for as many terms as possible.

Let's  review first how to prove that the number of order ideals of $P_{k+1,k+2}$, let's call it $p(k)$,  is the Catalan number $C_{k+1}$.
Let $i$ be the smallest empty label on the hypotenuse, implying
that $1, \dots, i-1$ are occupied, and `kicking out' all  vertices that are $\geq_P$ of the vertex labeled $i$,
leaving us with two connected components, triangles of sizes $i-2$ and $k-i$, with independent decisions regarding their order ideals. 
The `initial conditions' are $p(-1)=1$, $p(0)=1$,  and for $k \geq 1$, we have
$$
p(k)=\sum_{i=1}^{k+1} p(i-2) p(k-i)  .
\eqno(0)
$$

Now let's move-on to finding $s(n)$, i.e. the number of order ideals of $P_{2n+1,2n+3}$ without consecutive labels.

A diagram of the lattice $P_{2n+1,2n+3}$ (for $n=6$)  can be found in Figure 1(a) (see also Figure 3 (page 5) of [YQJZ], where the lattice is drawn such that the rank-zero vertices are at the bottom rather than on the diagonal).

%

\begin{figure}
\centering
\begin{subfigure}[t]{0.7\textwidth}
\includegraphics[width=\textwidth]{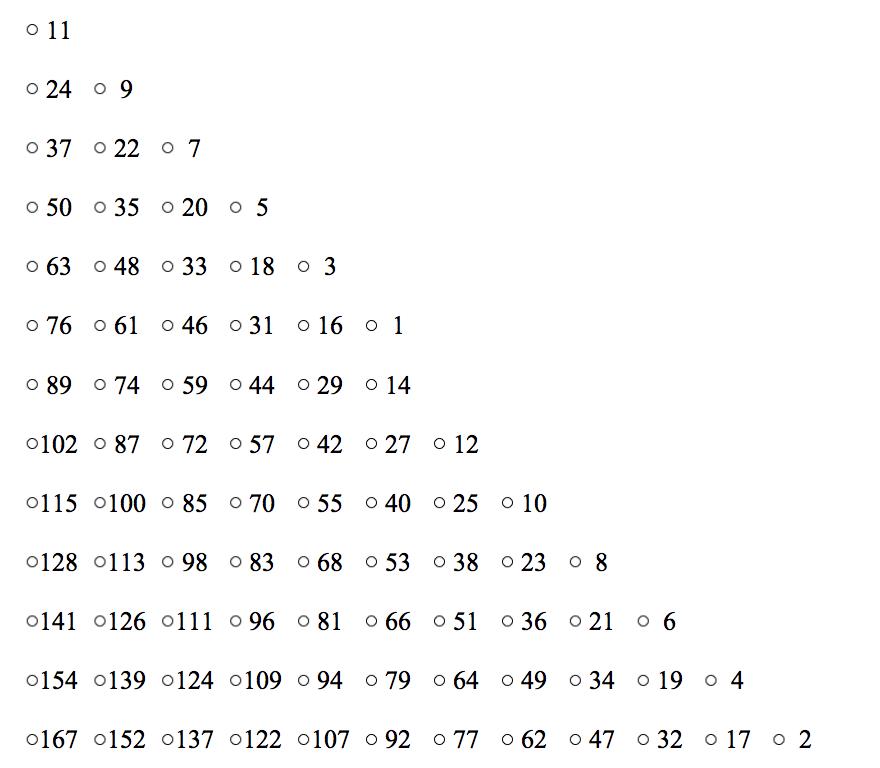}
\caption{The lattice $P_{13,15}$.}
\end{subfigure}
\begin{subfigure}[b]{0.7\textwidth}
\includegraphics[width=\textwidth]{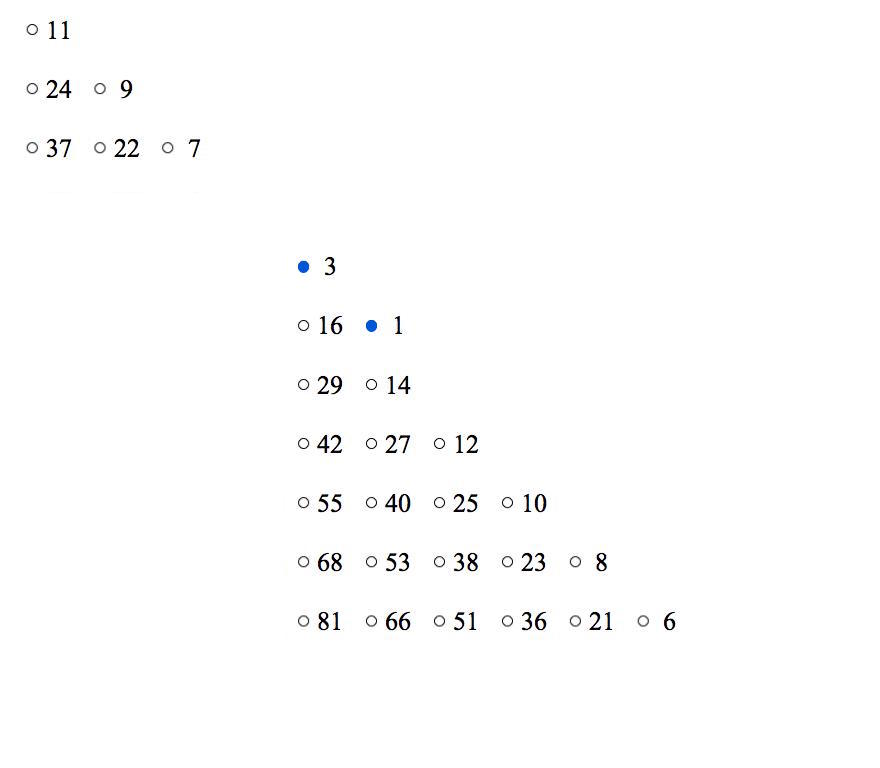}
\caption{A sub-lattice of $P_{13,15}$ which contains all order ideals of $P_{13,15}$ with smallest unoccupied odd label 5 and no consecutive labels. Note that this lattice is isomorphic to $EO(3,6)$ union the labels $1,3$.}
\end{subfigure}
\caption{}
\end{figure}

Inspired by the reasoning in [YQJZ],
let $2i-1$ ($1 \leq i \leq k$), be the smallest odd vertex (of rank $0$)
that is {\bf unoccupied}. This means that the vertices labeled $1, 3, \dots, 2i-3$ are {\bf occupied}.
This means that the  vertices with even labels, $2, \dots, 2i-2$ are {\bf unoccupied}, and since we are talking about {\it order ideals},
everything $\geq$ the odd vertex $2i-1$ and above the even vertices $2, \dots, 2i-2$ gets kicked out, and
for this scenario, we are left with counting order ideals of a smaller lattice, with two connected components,
that consists of an even-labeled component, a triangle-lattice whose rank zero level has size $n$, and whose labels are
$2i, 2i+2, \dots, 2i+2n-2$, and an odd-labeled component, a triangle whose rank zero level has  $n-i$ vertices, and whose  labels
are $2i+1, 2i+3, \dots, 2n-1$. In addition we have the definitely occupied vertices  $1, \dots, 2i-3$, but since they are definitely
occupied, they don't contribute anything to the count of order ideals.  

Figure 1(b)  depicts the case when labels $1$ and $3$ of $P_{13,15}$ are occupied and $5$ is empty.  All vertices $\geq$ $5,2,4$ cannot be part of the order ideal.

Let $EO(a,b)$ be a two-triangle lattice, consisting of a triangle  with $a$ rank-zero vertices
whose labels are $2, \dots , 2a$, and a triangle of length-side $b$ ($b>a$) whose labels are $1,3, \dots, 2b-1$. (See Figure 2(a) for a picture of $EO(7,9)$.)
Going back to the paragraph above, subtracting $2i-1$
from all labels, gives us a lattice isomorphic to $EO(n-i,n)$. Let $e(a,b)$ be the number of order ideals of
the lattice $EO(a,b)$ without consecutive labels. Then we have
$$
s(n)= \sum_{i=1}^{n+1} e(n-i,n)  .
\eqno(1)
$$

\begin{figure}
\centering
\begin{subfigure}[t]{0.7\textwidth}
\includegraphics[width=\textwidth]{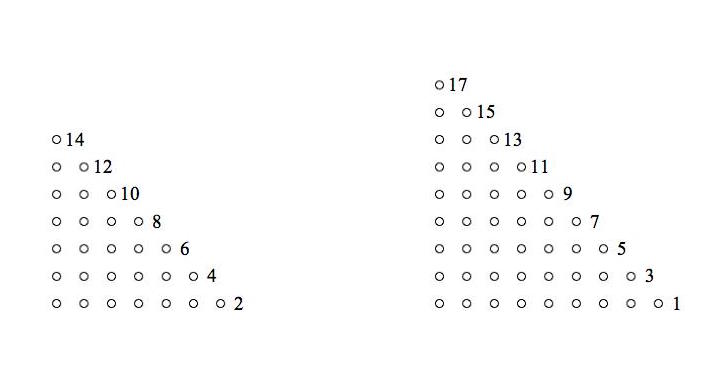}
\caption{The lattice $EO(7,9)$.}
\end{subfigure}
\begin{subfigure}[b]{0.7\textwidth}
\includegraphics[width=\textwidth]{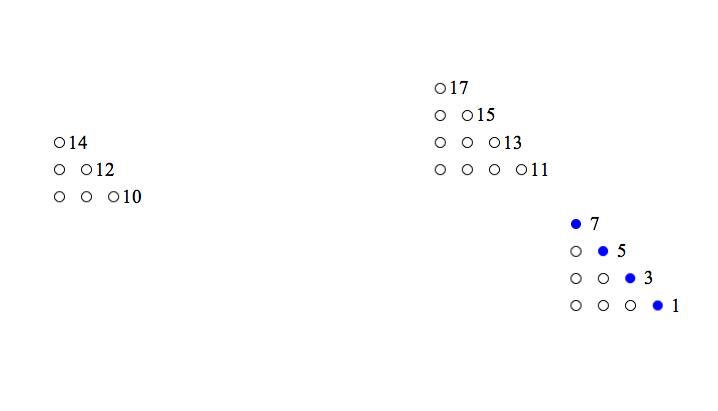}
\caption{A sub-lattice of $EO(7,9)$ which contains all order ideals of $EO(7,9)$ with smallest unoccupied odd label 9 and no consecutive labels.  Note that this lattice is isomorphic to $OE(3,4)$ union the triangular component containing the labels $1,3,5,7$.}
\end{subfigure}
\caption{}
\end{figure}

So if we would have an efficient `scheme' to compute $e(a,b)$, then we would be able to compute our sequence-of-desire $s(n)$.

For $a \leq b$, let  $OE(a,b)$ be $EO(b,a)$, and let $o(a,b)$ be the number of order ideals without consecutive labels of $OE(a,b)$.

By looking at the smallest unoccupied odd-labeled vertex, say $2i-1$  (see Figure 2(b)) we get, for $a \geq 1$:
$$
e(a,b)=\sum_{i=1}^{b+1} o(a+1-i,b-i) \, p(i-2)  ,
\eqno(2)
$$

and for $a\leq 0$, we have $e(a,b)=p(b)$.  
Similarly, for $a \geq 1$,
$$
o(a,b)=\sum_{i=1}^{a+1} e(a-i,b+1-i) \,       p(i-2)  ,
\eqno(3)
$$
and for $a \leq 0$, we have $o(a,b)=p(b)$.

The scheme consisting of equations $(0-3)$ enables a very fast computation of the sequence $s(i)$, for, say $i \leq 400$, confirming,
empirically for now, that $s(i)=4^i$. However this can be easily turned into a fully rigorous proof.
A {\it holonomic description} (see [Ze1], beautifully implemented by Christoph Koutschan in [K]) of both $e(a,b)$ and $o(a,b)$
can be readily guessed, and then, along with $p(k)=C_{k+1}$, the resulting identities $(1)-(3)$ are routinely verifiable
identities in the holonomic ansatz, that can be plugged into Koutschan's `holonomic calculator'. But since we know
{\it a priori} that $s(k)$ satisfies {\it some} such recurrence, and it is extremely unlikely that its order is very high,
confirming it for the first $400$ values consists a convincing {\it semi-rigorous proof}, that is easily rigorizable
(if [stupidly!]   desired).

{\bf Added in revised version:} Armin Straub found a far slicker, less computer-heavy, way,  to conclude this
experimental mathematics proof. See [Str2].

\subsection{Weight Enumerators}

But our main goal is to have $(2n+1,2n+3)$-analogs of the work in the article [Za] that dealt with $(n,n+1)$-core partitions with distinct parts.
In order to get data for the expectation, variance, and moments, we need an efficient way to generate as many terms of the
sequence of {\it Straub polynomials}, $S_n(q)$, defined by
$$
S_n(q) := \sum_{p} \,  q^{size(p)} ,
$$
where the sum ranges over all $(2n+1,2n+3)$-core partitions with distinct parts, $p$, and $size(p)$ is the sum of the
entries of $p$ (i.e. the number of boxes in its Young Diagram).

The Maple package {\tt core.txt} that accompanied [Za], and is also accompanying this article, uses
Dyck paths, and was able to find the first nine Straub polynomials, $S_n(q)$, $1 \leq n \leq 9$. 
It is based on an extension of the method described in [EZ], but keeping track of the fact that
cells with adjacent labels are not allowed. So one has to put up with much more general
families of paths, that are also parametrized by a set of `forbidden labels'. This causes an exponential
expansion of memory and time.

The approach that we take in this article, that easily produced the first $21$ Straub polynomials,
is a weighted analog of the above naive-enumeration scheme, and goes via order ideals.

For an order ideal of $P_{m,n}$ let its {\it weight} be
$$
q^{SumOfLabels} t^{NumberOfVertices}  .
$$

Let $Q(n)$ be the set of order ideals of $P_{2n+1,2n+3}$ without neighboring labels (i.e. if $a \in I$ then
both $a-1$ and $a+1$ are not in $I$). Let's define the two-variable polynomials
$$
A_n(q,t):=\sum_{I \in Q(n)} q^{SumOfLabels(I)} t^{NumberOfVertices(I)}  .
$$
Define the `umbra' (linear functional on polynomials of $t$) by
$$
U(t^k):=q^{-k(k-1)/2} ,
$$
and extended linearly.
As shown by Anderson, once $A_n(q,t)$ are known, we get $S_n(q)$ by
$$
S_n(q)= U(A_n(q,t)) ,
$$
in other words, to get $S_n(q)$ replace any power,  $t^k$, that appears in $A_n(q,t)$, by $q^{-k(k-1)/2}$.

It remains to find an efficient scheme for `cranking out' as many terms of $A_n(q,t)$ that our
computer would be willing to compute.

We first need a weighted analog of Equation $(0)$, i.e. the weight-enumerator of $P_{k+1,k+2}$, but we need the
extra generality where  (still with the smallest label being $1$), for {\it any} positive integers $c$ and $h$,
in the vertical direction it is going down by $c$, and in the horizontal direction it going down by $c+h$
(drawing the lattice so that the highest label, $1+(c+h)(k-1)$ is at the origin, and the vertex labeled $1$
is situated at the point $(k-1,0)$, and the vertex labeled $1+(k-1)h$ is situated at the
point $(0,k-1)$. Note that the original $P_{k+1,k+2}$ corresponds to $c=k+1$ and $h=1$.

Let's call this generalized weight-enumerator $P^{(c,h)}_k (q,t)$. It is readily seen that the weighted analog of Eq. $(0)$ is
$$
P^{(c,h)}_k (q,t) \, = \,
\sum_{i=1}^{k+1} \,
t^{i-1} \, \cdot \, q^{(i-1)+(i-1)(i-2)h/2} \, \cdot \,  P^{(c,h)}_{i-2} (q, q^{c+h} t) \, \cdot \, P^{(c,h)}_{k-i} (q, q^{ih} t) ,
\eqno(0w)
$$
with the initial conditions $P_{-1}=1, P_{0}=1$.

Let $E^{(c)}_{x,y} (q,t) $ be the weight-enumerator of the lattice $EO(x,y)$
with  horizontal spacing $c$ and vertical spacing $c+2$.  Then the analog of Eq. $(1)$ is
$$
A_n(q,t) \, = \,
\sum_{i=1}^{n+1} t^{i-1} q^{(i-1)^2} \, \cdot \,        E^{(2n+1)}_{n-i,n} (q,q^{2i-1} t)  .
\eqno(1w)
$$

Let $O^{(c)}_{x,y} (q,t) $ be the weight-enumerator of the lattice $OE(x,y)$,
with  horizontal spacing $c$ and vertical spacing $c+2$.  Then the analog of Eq. $(2)$ 
can be seen to be
$$
E^{(c)}_{x,y}(q,t) \, = \, \sum_{i=1}^{y+1} \,
t^{i-1}\, \cdot \, q^{(i-1)^2} \, \cdot \,  O^{(c)}_{x-i+1,y-i} (q, q^{2i-1} t) \, \cdot \,  P^{(c,2)}_{i-2} (q,q^{c+2}t)
 ,
\eqno(2w) 
$$
with the {\it initial condition} 
$E^{(c)}_{x,y}(q,t)=P^{(c,2)}_y(q,t)$ when $x \leq 0$.

Finally, the weighted analog of Eq. $(3)$ is
$$
O^{(c)}_{x,y}(q,t) \, = \, \sum_{i=1}^{x+1} \,
t^{i-1}\, q^{(i-1)^2} \, \cdot \,  E^{(c)}_{x-i,y-i+1} (q, q^{2i-1} t) \,  \cdot \, P^{(c,2)}_{i-2} (q,q^{c+2}t) ,
\eqno(3w)
$$
with the {\it initial condition}  $O^{(c)}_{x,y}(q,t)=P^{(c,2)}_y(q,qt)$ when $x \leq 0$.

\subsection{The first $21$  Straub polynomials}

Using the above scheme, one gets that
$$
S_1(q) \, = \, {q}^{4}+{q}^{2}+q+1  ,
$$
$$
S_2(q) \, = \,  {q}^{21}+{q}^{16}+2\,{q}^{12}+{q}^{9}+{q}^{8}+{q}^{7}+{q}^{6}+{q}^{5}+2\,{q}^{4}+2\,{q}^{3}+{q}^{2}+q+1
$$
$$
S_3(q) \, = \,  
{q}^{65}+{q}^{56}+{q}^{48}+{q}^{47}+{q}^{41}+{q}^{39}+{q}^{37}+2\,{q}^{35}+{q}^{32}+{q}^{30}+2\,{q}^{29}+{q}^{28}+{q}^{26}+3\,{q}^{24}+{q}^{23}+{q}^{22}
$$
$$
+{q}^{21}+{q}^{20}+2\,{q}^{19}+2\,{q}^{18}+3\,{q}^{17}+{q}^{16}+{q}^{15}+2\,{q}^{14}+2\,{q}^{13}+2\,{q}^{12}+3\,{q}^{11}+{q}^{10}+3\,{q}^{9}+3\,{q}^{8}
$$
$$
+3\,{q}^{7}+4\,{q}^{6}+3\,{q}^{5}+2\,{q}^{4}+2\,{q}^{3}+{q}^{2}+q+1  ,
$$
$$
S_4(q) \, = \,  
{q}^{155}+{q}^{141}+{q}^{128}+{q}^{125}+{q}^{116}+{q}^{112}+2\,{q}^{105}+{q}^{103}+{q}^{100}+2\,{q}^{95}+{q}^{93}+{q}^{91}+2\,{q}^{89}+{q}^{85}+{q}^{84}
$$
$$
+{q}^{83}+2\,{q}^{82}+{q}^{80}+{q}^{79}+{q}^{78}+{q}^{76}+{q}^{74}+{q}^{73}+{q}^{72}+2\,{q}^{71}+2\,{q}^{70}+{q}^{69}+2\,{q}^{68}+{q}^{67}+{q}^{65}+{q}^{64}
$$
$$
+{q}^{63}+5\,{q}^{61}+{q}^{60}+2\,{q}^{59}+3\,{q}^{57}+{q}^{56}+3\,{q}^{55}+4\,{q}^{53}+2\,{q}^{52}+2\,{q}^{51}+2\,{q}^{50}+{q}^{49}+2\,{q}^{48}+3\,{q}^{47}
$$
$$
+2\,{q}^{46}+3\,{q}^{45}+4\,{q}^{44}+2\,{q}^{43}+{q}^{42}+5\,{q}^{40}+3\,{q}^{39}+4\,{q}^{38}+5\,{q}^{37}+2\,{q}^{36}+3\,{q}^{35}+{q}^{34}+4\,{q}^{33}
$$
$$
+6\,{q}^{32}+5\,{q}^{31}+3\,{q}^{30}+4\,{q}^{29}+3\,{q}^{28}+5\,{q}^{27}+4\,{q}^{26}+7\,{q}^{25}+5\,{q}^{24}+6\,{q}^{23}+3\,{q}^{22}+4\,{q}^{21}+5\,{q}^{20}
$$
$$
+5\,{q}^{19}+4\,{q}^{18}+5\,{q}^{17}+6\,{q}^{16}+5\,{q}^{15}+4\,{q}^{14}+7\,{q}^{13}+6\,{q}^{12}+7\,{q}^{11}+7\,{q}^{10}+6\,{q}^{9}+6\,{q}^{8}+5\,{q}^{7}
$$
$$
+4\,{q}^{6}+3\,{q}^{5}+2\,{q}^{4}+2\,{q}^{3}+{q}^{2}+q+1 .
$$

For the Straub polynomials $S_n(q)$ for $5 \leq n \leq 21$, see the webpage  \hfill\break
{\tt http://www.math.rutgers.edu/\~{}zeilberg/tokhniot/oArmin3.txt} , or use procedure {\tt ASpc(n,q)} in the Maple package {\tt Armin.txt}
mentioned above.


Unlike the case of $(s,s+1)$-core partitions, whose number happened to be $F_{s+1}$, and  the
explicit expressions for the expectation, variance, and higher moments involved expressions in
$F_s,F_{s+1}$ and $s$, the present case of $(2n+1,2n+3)$-core partitions into distinct parts,
gives, surprisingly, `nicer' results.
This is because, as conjectured in [Str1] and first proved in [YQLZ] (and reproved above), the actual enumeration
is as simple as can be, namely $4^n$. Hence it is not surprising that the expectation, variance, and higher
moments are {\it polynomials} in $n$.

To get expressions for the moments we used the empirical-yet-rigorizable approach of [Ze2] and [Ze3], as follows.

Using the first $21$ Straub polynomials, we get the sequence  
of numerical averages $S_n'(1)/4^n$, $1 \leq n \leq 21$,
and `fit it' to a polynomial of degree $3$ (in fact four terms suffice!),
we get the expression for the expectation, let's call it $\mu(n)$,
stated in Theorem 1 above.

Using the sequence 
$$
\frac{(q \frac{d}{dq})^2 S_n(q)|_{q=1}}{4^n} -\mu(n)^2  ,
$$
for $1 \leq n \leq 7$, and `fitting' it  with a polynomial of degree $6$, we get an explicit expression for the variance, thereby getting Theorem 2.
The conjectured polynomial expression agrees all the way to $n=21$.

The third-through the seventh moments are derived similarly, where the $i$-th moment (about the mean, but also the straight moment)
turns out to be a polynomial of degree $3i$ in $n$.

Let us comment that we strongly believe that
all the results here can be, {\it a posteriori}, justified rigorously. The complicated functional recurrences
for the Straub polynomials (before the ``umbral application'') entail, after Taylor expansions about $q=1$, extremely
complicated recurrence relations for the (pre-) moments, whose details do not concern us, since we know that
their truth follows by induction. The reason that we are not completely sure about this is that we don't have
a formal proof that ``polynomiality'' is preserved under the umbral transform.
Granting this, each such identity is a {\it polynomial identity}, and hence its truth follows
from plugging-in sufficiently many special cases. But that's how we got them in the first place. {\bf QED}!

\subsection{Encore: A one-line (almost) proof of Straub's Ex-Conjecture about the Maximal Size of a $(2n+1,2n+3)$ core partition into distinct parts}

In [YQLZ], the authors used quite a bit of {\it human ingenuity} to prove Armin Straub's conjecture (posed in [Str1])
that the maximal size of a $(2n+1,2n+3)$-core partition into distinct parts is given by the degree-$4$ polynomial
$\frac{1}{24}\, \left( 5\,n+11 \right) n \left( n+2 \right)  \left( n+1 \right)$.

We strongly believe that one can deduce from general, {\it a priori}, {\it hand-waving} (yet fully rigorous) considerations that this quantity is
{\it some} polynomial of degree $\leq 5$. Hence it is enough to check it for $1 \leq n \leq 6$. But
this quantity is exactly the {\bf degree} of the Straub polynomial $S_n(q)$. We verified it, in fact, all the  way to
$n=21$, so Theorem $0'$ is re-proved (modulo our belief) (with a {\it vengeance}!).

\section*{Acknowledgment}

We wish to thank Tewodros Amdeberhan for introducing us to this interesting subject. 
We also wish to thank Armin Straub for carefully reading an earlier draft, and making many insightful comments.

\section*{References}

[Am] Tewodros Amdeberhan, {\it Theorems, problems and conjectures},\\ {\tt https://arxiv.org/abs/1207.4045}  .

[AmL] Tewodros Amdeberhan and  Emily Sergel Leven , {\it Multi-cores, posets, and lattice paths}, \hfill\break
{\tt https://arxiv.org/abs/1406.2250}.
Also published in Adv. Appl. Math. {\bf 71}(2015), 1-13.

[An] Jaclyn Anderson, {\it Partitions which are simultaneously $t_1$  and $t_2$ -core}, Discrete Math. {\bf 248}(2002), 237-243.

[AHJ] Drew Armstrong, Christopher R.H. Hanusa, and B. Jones, {\it Results and conjectures on simultaneous core partitions}, 
{\tt https://arxiv.org/abs/1308.0572}  .
Also published in European J. Combin. {\bf 41} (2014), 205-220.

[EZ] Shalosh B. Ekhad and Doron Zeilberger,
{\it Explicit Expressions for the Variance and Higher Moments of the Size of a Simultaneous Core Partition and its Limiting Distribution},
The Personal Journal of Shalosh B. Ekhad and Doron Zeilberger, \hfill\break
{\tt http://www.math.rutgers.edu/~zeilberg/pj.html}  . \hfill\break
Direct url: \\{\tt http://www.math.rutgers.edu/~zeilberg/mamarim/mamarimhtml/stcore.html}. 
Also in: {\tt https://arxiv.org/abs/1508.07637}  .

[J] Paul Johnson, {\it Lattice points and simultaneous core partitions}, \hfill\break
{\tt http://arxiv.org/abs/1502.07934}  , 27 Feb 2015.

[K] Christoph Koutschan, {\it Advanced applications of the holonomic systems approach},
Research Institute for Symbolic Computation (RISC), Johannes Kepler University, Linz, Austria, 2009.  \hfill\break
{\tt http://www.koutschan.de/publ/Koutschan09/thesisKoutschan.pdf}. \hfill\break
Software packages available from: \hfill\break
{\tt http://www.risc.jku.at/research/combinat/software/HolonomicFunctions/}.

[StaZ] Richard P. Stanley and Fabrizio Zanello, {\it The Catalan case of Armstrong's conjecture on simultaneous core partitions},  \hfill\break
{\tt http://arxiv.org/abs/1312.4352}. Also published in:
SIAM J. Discrete Math. {\bf 29}(2015) , 658-666.

[Str1] Armin Straub, {\it Core partitions into distinct parts and an analog of Euler's theorem}, \hfill\break
{\tt https://arxiv.org/abs/1601.07161}.
Also published in: European J. of Combinatorics {\bf 57} (2016), 40-49.

[Str2] Armin Straub, {\it Email message to the authors}, \hfill\break
{\tt  http://www.math.rutgers.edu/\~{}zeilberg/mamarim/mamarimhtml/arminFeedback.pdf}.

[TW] Marko Thiel and Nathan Williams, {\it Strange Expectations},\\ {\tt https://arxiv.org/abs/1508.05293}, 21 Aug. 2015.

[Wan] Victor Y. Wang, {\it Simultaneous core partitions: parameterizations and sums}, \hfill\break
{\tt http://arxiv.org/abs/1507.04290}. Also published in:
Electron. J. Combin.  {\bf 23}(2016), Paper 1.4, 34 pp.

[X] Huan Xiong,  {\it  Core partitions with distinct parts}, \hfill\break
{\tt https://arxiv.org/abs/1508.07918}.

[YQJZ] Sherry H.F. Yan, Guizhi Qin, Zemin Jin, Robin D.P. Zhou, {\it On $(2k+1,2k+3)$-core partitions with distinct parts}, 
\hfill\break
{\tt https://arxiv.org/abs/1604.03729}.

[Za] Anthony Zaleski, {\it Explicit expressions for the moments of the size of an $(s,s+1)$-core partition with distinct parts},
{\tt https://arxiv.org/abs/1608.02262}. Also appeared in Advances in Applied Mathematics {\bf 84} (2017), 1-7.

[Ze1] Doron Zeilberger, {\it A Holonomic Systems Approach To Special Functions},
J. Computational and Applied Math {\bf 32} (1990), 321-368, \hfill\break
{\tt http://www.math.rutgers.edu/\~{}zeilberg/mamarim/mamarimPDF/holonomic.pdf}.

[Ze2] Doron Zeilberger, {\it The Automatic Central Limit Theorems Generator (and Much More!)},
in: ``Advances in Combinatorial Mathematics: Proceedings of the Waterloo Workshop in 
Computer Algebra 2008 in honor of Georgy P. Egorychev'', chapter 8, pp. 165-174, (I.Kotsireas, E.Zima, eds. Springer Verlag, 2009),
\hfill\break
{\tt http://www.math.rutgers.edu/\~{}zeilberg/mamarim/mamarimhtml/georgy.html}.

[Ze3]  Doron Zeilberger, {\it HISTABRUT: A Maple Package for Symbol-Crunching in Probability theory},
the Personal Journal of Shalosh B. Ekhad and Doron Zeilberger, posted Aug. 25, 2010,
\hfill\break
{\tt http://www.math.rutgers.edu/\~{}zeilberg/mamarim/mamarimhtml/histabrut.html}.

[Ze4]  Doron Zeilberger, {\it Symbolic Moment Calculus I.: Foundations and Permutation Pattern Statistics},
Annals of Combinatorics {\bf 8} (2004), 369-378. Available from:
\hfill\break
{\tt http://www.math.rutgers.edu/\~{}zeilberg/mamarim/mamarimhtml/smcI.html}.

\bigskip
\hrule
\bigskip
Anthony Zaleski, Department of Mathematics, Rutgers University (New Brunswick), Hill Center-Busch Campus, 110 Frelinghuysen
Rd., Piscataway, NJ 08854-8019, USA. \hfill\break
Email: {\tt az202 at math dot rutgers dot edu}    .

Doron Zeilberger, Department of Mathematics, Rutgers University (New Brunswick), Hill Center-Busch Campus, 110 Frelinghuysen
Rd., Piscataway, NJ 08854-8019, USA. \hfill\break
Email: {\tt DoronZeil at gmail dot com}    .

\end{document}